\documentclass[12pt]{article}
\usepackage{amsmath, amssymb,amsthm}
\usepackage{graphicx}
\usepackage{amsthm}
\usepackage{pstricks-add}

\usepackage{enumerate}

\newcommand*{\F}{\mathcal{F}}
\newcommand*{\B}{\mathcal{B}}

\newcommand*{\A}{\mathcal{A}}

\newtheorem{theorem}{Theorem}
\newtheorem{lemma}{Lemma}
\newtheorem{Fact}{Fact}

\begin{document}

\title{Short Proof of Erd\H os Conjecture for Triple Systems}

\author{ Peter Frankl\\R\'enyi Institute\\Budapest, Hungary\\{\tt peter.frankl@gmail.com}
\and Vojtech R\"odl\thanks{Research supported  by NSF grant DMS 1301698}
\\Emory University\\Atlanta, GA\\{\tt rodl@mathcs.emory.edu}
\and Andrzej Ruci\'nski
\thanks{Research supported by the Polish NSC grant 2014/15/B/ST1/01688.
Part of research performed at Emory University, Atlanta.}
\\A. Mickiewicz University\\Pozna\'n, Poland\\{\tt rucinski@amu.edu.pl}
}

\date{\today}

\maketitle

\let\thefootnote\relax
\let\thefootnote\relax

\begin{abstract}
In 1965 Erd\H
os  conjectured that for all $k\ge2$, $s\ge1$ and $n\ge k(s+1)$, an $n$-vertex $k$-uniform
hypergraph $\F$ with $\nu(\F)=s$ cannot have more than
\newline
$\max\{\binom{sk+k-1}k,\;\binom nk-\binom{n-s}k\}$ edges. It took almost
fifty years to prove it for  triple systems. In 2012  we proved the conjecture for all $s$ and all $n\ge4(s+1)$. Then {\L}uczak and Mieczkowska
(2013) proved the conjecture for sufficiently large $s$ and all $n$. Soon after, Frankl proved it for all $s$. Here we present a simpler version of that  proof which yields Erd\H os's conjecture for $s\ge33$. Our motivation is to lay down foundations for a possible proof in the much harder  case $k=4$, at least for
large $s$.
\end{abstract}

\section{Introduction}\label{intro}

Let $\nu(\F)$ denote the size of a largest matching in a $k$-uniform hypergraph $\F$. In 1965 Erd\H
os \cite{E} conjectured that for all $k\ge2$, $s\ge1$ and $n\ge k(s+1)$, an $n$-vertex $k$-uniform
graph $\F$ with $\nu(\F)=s$ cannot have more than
$$\max\left\{\binom{sk+k-1}k,\;\binom
nk-\binom{n-s}k\right\}$$ edges and proved it for $n> n_0(k,s)$. So far, the best general
upper bound,  $n_0(k,s)\le (2s+1)k-s$, is due to Frankl \cite{F-JCTA}. While the conjecture was proved
for $k=2$ (i.e., for graphs) by Erd\H os and Gallai already in 1959 \cite{EG}, the progress for $k\ge3$ has been much slower.

 It took almost
fifty years to prove it for  triple systems ($k=3$). First, in \cite{FRR} we proved the conjecture for all $s$ and all $n\ge4(s+1)$. Then {\L}uczak and Mieczkowska \cite{LM} proved the conjecture for sufficiently large $s$ and all $n$. Soon after, Frankl \cite{F-arxiv} proved it for all $s$,
building upon some ideas from \cite{LM}. Here we present a streamlined version of the  proof in \cite{F-arxiv} which
is shorter and simpler and yields the conjecture for $s\ge33$. Our motivation for writing this paper
is also rooted at the belief that a proof along the same lines could  eventually work in the case $k=4$, at least for
large $s$. At one place in the proof we rely on the above mentioned result from \cite{FRR}, which for $k\ge4$ could be replaced by the general result from \cite{F-JCTA} (also mentioned above).

We will call a 3-uniform hypergraph a \emph{3-graph}.
Let
$$m(n,s)=\max\{|\F|: |V(\F)|=n,\; \nu(\F)=s\}.$$
Further, let
$$\A:=\A(n,s)=K^3_{3s+2}\cup (n-3s+2)K_1$$ be the complete 3-graph  on $3s+2$ vertices, augmented by $n-3s+2$ isolated
vertices, and let
$$\B:=\B(n,k)=K^3_{n}-K^3_{n-s},$$ that is, $\B$ is the complete 3-graph on $n$ vertices from
which a complete 3-graph on $n-s$ vertices has been removed.
Finally,
let $$a(s)=|\A(n,s)|=\binom{3s+2}3\quad\mbox{ and }\quad b(n,s)=|\B(n,s)|=\binom n3-\binom{n-s}3$$
and
$$M(n,s)=\max\left\{a(s),\;b(n,s)\right\}.$$

 Clearly, $m(n,s)\ge M(n,s)$, and for $k=3$ the
Erd\H os Conjecture states that $m(n,s)=M(n,s)$ for all $s$ and all $n$.  Note that the conjecture is trivially valid for $n\le 3s+2$, since then $m(n,s)\le a(s)$. Therefore, it suffices to consider only the case $n\ge3s+3$. Here we prove the
following result.

\begin{theorem}\label{main}
For all $s\ge33$ and $n\ge 3s+3$, $m(n,s)=M(n,s)$.
\end{theorem}

\section{Preparations for the proof}
Given a linear order $\le$ on the vertex set $V(\F)$, we  define a partial order on the triples of
vertices as follows.
 For two sets $A,B\in\binom{[n]}3$ we write $A\prec B$ if $A=\{a_1\le a_2\le a_3\}$,
$B=\{b_1\le b_2\le b_3\}$, and $a_i\le b_i$ for all $i=1,2,3$. We say that  $\F$ is \emph{stable
(or shifted)} if whenever $A\prec B$ and $B\in\mathcal F$, then $A\in\mathcal F$. If a 3-graph $\F$
is not stable, there exists an edge $B$ of $\F$ and a non-edge $A$ with $A\prec B$. Then we might
swap $A$ and $B$ in $\F$. We call such an operation \emph{a shift}.

Let $sh(\F)$ be a stable
3-graph obtained from $\F$ by a series of shifts. Note that $|sh(\F)|=|\F|$. It is an easy
exercise (see, e.g., \cite[Lemma 3]{LM}) to show that $\nu(sh(\F))\le \nu(\F)$. This means
that it is sufficient to prove the Erd\H os Conjecture for stable 3-graphs only.
For further reference, note also that if
$$\nu(\F)=s\quad\mbox{and}\quad|\F|=m(n,s),\quad\mbox{then}\quad\nu(sh(\F))=\nu(\F).$$

When comparing the two quantities defining $M(n,s)$, it is apparent that for smaller~$n$ we have
$a(s)>b(n,s)$, while for larger $n$ the opposite holds. Indeed, $b(n,s)$ is an increasing function of $n$, while $a(s)$ is a constant.
For every $s$, we define
$$n_1(s)=\min\{n: a(s)\le b(n,s)\}.$$
It follows that $M(n,s)=a(s)$ for $n\le n_1(s)-1$, while $M(n,s)=b(n,s)$ for $n\ge n_1(s)$.
The importance of the parameter $n_1(s)$ is facilitated by the  fact, observed already in
\cite[Fact 5.3]{F-arxiv}, that if the Erd\H os conjecture fails for some $s$ and $n$, then it must fail for that $s$ and $n\in\{n_1(s)-1,n_1(s)\}$. We provide the proof for completeness. Given a 3-graph $\F$ and a vertex
$v\in V(\F)$,  let
$$\F(\bar v)=\{F\in\F:\;v\not\in F\}\quad\mbox{ and }\quad \F(v)=\{F\setminus\{v\}:\;v\in
F\in\F\}.$$ Note that $\F(\bar v)$ is a 3-graph, while $\F(v)$ is a graph, both on the same vertex set \newline $V(\F)\setminus\{v\}$.

\begin{Fact}\label{F1}
If for some $s$,  $m(n_1(s)-1,s)=M(n_1(s)-1,s)$ and $m(n_1(s),s)=M(n_1(s),s)$, then
$m(n,s)=M(n,s)$ for all $n\ge3s+3$.
\end{Fact}
\proof Assume that
$$m(n_1(s)-1,s)=M(n_1(s)-1,s)=a(s)\quad\mbox{and}\quad m(n_1(s),s)=M(n_1(s),s)=b(n_1(s),s).$$ Since $a(s)$
is independent of $n$, it follows that $m(n,s)=a(s)$ for all $n\le n_1(s)-1$. For $n\ge n_1(s)$ we
use induction on $n$. Assume $m(n-1,s)=b(n-1,s)$, $n\ge n_1(s)+1$, and let $\F$
be a stable 3-graph on vertex set $[n]=\{1,\dots,n\}$ (ordered by natural ordering) and with
$\nu(\F)=s$. Consider the 3-graph $\F(\bar n)$ and the graph $\F(n)$. Clearly, \newline $|\F|=|\F(\bar
n)|+|\F(n)|$ and  $\nu(\F(\bar n))\le \nu(\F)=s$, so, by assumption,
$$|\F(\bar n|\le
m(n-1,s)=b(n-1,s)=\binom{n-1}3-\binom{n-1-s}3.$$

By stability of $\F$, we also have $\nu(\F( n))\le s$. Indeed, if there was a matching
$e_1,\dots,e_{s+1}$ in $\F(n)$, then, since $n-1\ge3s+3$, the triples $e_i\cup\{v_i\}$,
$i=1,\dots,s+1$, where $\{v_1,\dots,v_{s+1}\}\subset [n-1]\setminus\bigcup_{i=1}^{s+1}e_i$, would
form a matching of size $s+1$ in $\F$, a contradiction.  Thus, by the result of Erd\H os and Gallai
from \cite{EG} quoted earlier,
\begin{equation}\label{graphs}
|\F(n)|\le \max\left\{\binom{2s+1}2,\;\binom {n-1}2-\binom{n-1-s}2\right\}=
\binom{n-1}2-\binom{n-1-s}2
\end{equation}
for $n\ge3s$.
Altogether,
$$|\F|\le\binom{n-1}3-\binom{n-1-s}3+\binom{n-1}2-\binom{n-1-s}2=\binom n3-\binom{n-s}2=b(n,s). \qed$$

\medskip

The value of $n_1(s)$ was asymptotically determined in \cite{LM} and \cite{FRR}. Here we need
estimates valid also for small $s$.

\begin{Fact}\label{n0s} For all $s$ and $n$,
  $${\rm{(i)}}\quad n_1(s)\le 3.5s+3,\quad\rm{(ii)}\quad n_1(s)\ge 3.4s+1,\quad\rm{(iii)}\quad n_1(s)-n_1(s-1)\ge2.$$
\end{Fact}

\proof All parts follow from an exact formula for $n_1(s)$ which we derive first. We look for the
smallest integral solution (in $n$) to the inequality

$$\binom{3s+2}3\le\binom n3-\binom{n-s}3=\binom s3+\binom s2(n-s)+s\binom{n-s}2$$
which, after substituting $m=n-s$, becomes
$$3m^2+3(s-2)m-(26s^2+30s+4)\ge0.$$
Solving this quadratic inequality and setting $g(s)=321s^2+324s+84$, we derive that
\begin{equation}\label{n0s}
n_1(s)=1+\left\lceil\frac12s+\frac16\sqrt{g(s)}\right\rceil.
\end{equation}

To show (i), observe that, in view of (\ref{n0s}), it is now equivalent to
$$n_1(s)-1=\left\lceil\frac12s+\frac16\sqrt{g(s)}\right\rceil\le\lfloor{3.5s+2}\rfloor,$$
which, in turn, is equivalent to
$$\frac12s+\frac16\sqrt{g(s)}\le\lfloor{3.5s+2}\rfloor.$$
Since $3.5s+3/2\le\lfloor{3.5s+2}\rfloor$,
 we thus get a stronger inequality
$$g(s)\le36(3s+3/2)^2=324s^2+324s+81$$
which is true for all $s\ge1$.

Part (ii) is even easier and we leave it to the reader. For part (iii), rewriting $n_1(s-1)\le
n_1(s)-2$,
 substituting formula (\ref{n0s}), and dropping the ceilings on both sides, we obtain a stronger inequality
$$1+\frac12(s-1)+\sqrt{g(s-1)}\le 1+\frac12s+\sqrt{g(s)}-2,
$$
equivalent to
$$\frac{642s+3}{6(\sqrt{g(s-1)}+\sqrt{g(s)})}\ge\frac32.
$$
By bounding the L-H-S from below by
$$\frac{642s+3}{12\sqrt{g(s)}}$$
we finally obtain an yet stronger inequality
$$(214s+1)^2\ge36g(s),$$
valid for all $s\ge1$. \qed

\bigskip

We say that a stable 3-graph has property ONE if $\nu(\F(\bar1))=\nu(\F)$, that is, if there is a largest matching not covering the smallest vertex.
Let
$$m_{ONE}(n,s)=\max\{|\F|: |V(\F)|=n, \nu(\F)=s, \mbox{ and $\F$ has ONE }\}.$$
Clearly, $$m_{ONE}(n,s)\le m(n,s).$$ Note, however, that while  $\A$ has ONE, $\B$ does not. This
means that the inequality $m_{ONE}(n,s)\ge M(n,s)$ might not be true in general. On the other hand,
we are going to prove that the reverse inequality is true.

\begin{lemma}\label{L1}
$m_{ONE}(n,s)\le M(n,s)$ for all $s\ge25$ and all $n\ge3s+3$.
\end{lemma}

Lemma \ref{L1} is the main ingredient of our proof of Theorem \ref{main}.
The other ingredient is the following lemma.
\begin{lemma}\label{L2}
If, for some $s_0$, $m_{ONE}(n,s)\le M(n,s)$ for all $s\ge s_0$ and all $n\ge3s+3$, then $m(n,s)=M(n,s)$ for all
\begin{equation}\label{1}
s\ge\tfrac54(s_0+1)
\end{equation}
and all $n\ge3s+3$.
\end{lemma}

{\bf Proof of Theorem \ref{main}:\;}
By Lemma \ref{L1}, the assumption of Lemma \ref{L2}  is satisfied with $s_0=25$. Then Theorem \ref{main} follows from Lemma \ref{L2} with $s_0=25$, as the R-H-S of (\ref{1}) equals $\tfrac{65}2$. \qed

\medskip

 The proof of Lemma \ref{L2} is given below, while the proof of Lemma \ref{L1} is deferred to Section \ref{P1}.

\subsection{Proof of Lemma \ref{L2}}\label{P2}

The proof is based on a fact similar to Fact \ref{F1}.
\begin{Fact}\label{R}
For every $s$, if $m_{ONE}(n_1(s)-1,s)\le a(s)$ and $m_{ONE}(n_1(s),s)\le b(n_1(s),s)$, then
$m_{ONE}(n,s)\le M(n,s)$ for all $n\ge 3s+3$.
\end{Fact}

\proof For $n\le n_1(s)-2$, consider a stable 3-graph $\F$  on $n$ vertices, with property ONE,
with $|\F|=m_{ONE}(n,s)$, and with $\nu(\F)=s$. By adding to it vertices $n+1,\dots, n_1(s)-1$, we
obtain a 3-graph $\F'$ which is still stable, has property ONE and $\nu(\F)=s$. Thus
$$m_{ONE}(n,s)=|\F|=|\F'|\le m_{ONE}(n_1(s)-1,s)=a(s)=M(n,s).$$
For $n\ge n_1(s)+1$, we apply induction on $n$. Assume that $m_{ONE}(n-1,s)\le M(n-1,s)$ and let
$\F$, with $|\F|=m_{ONE}(n,s)$, be a stable 3-graph which has property ONE. To show that
$m_{ONE}(n,s)\le M(n,s)$, we proceed
 as in  the proof of Fact \ref{F1}.
The only novelty is to observe that if $\F$  has property ONE, then the same is true for $\F(\bar
n)$. This follows, since, due to stability of $\F$ and the fact that $n\ge3s+3$, there is a
matching of size $s$ in $\F$, avoiding vertices 1 and $n$. Hence, $|\F(\bar n)|\le
m_{ONE}(n-1,s)\le M(n-1,s)\le b(n-1,s)$. As we also have (\ref{graphs}), the conclusion follows.
\qed

\bigskip

{\bf Proof of Lemma \ref{L2}:} Suppose that for some $s\ge\tfrac54(s_0+1)$ and $n\ge 3s+3$, we have
$m(n,s)>M(n,s)$. By Fact \ref{F1}, it means that there is $n\in\{n_1(s)-1,n_1(s)\}$ and a 3-graph
$\F$ with $V(\F)=[n]$, $\nu(\F)=s$, and $|\F|=m(n,s)>M(n,s)$. Let $\F'=sh(\F)$. Recall that
$|\F'|=|\F|$ and $\nu(\F')=\nu(\F)=s$. Since $s\ge s_0$, $\F'$ cannot have property ONE, since
otherwise we would arrive at a contradiction with the assumption of Lemma 2.

For $q\ge1$, denote by  $\F_q'$ the induced sub-3-graph of $\F'$ obtained by removing the vertices
$1,2,\dots,q$ and all edges adjacent to them. Observe that $\F'_q$ is stable and $\nu(\F'_q)\ge
s-q$. Observe also that $\nu(\F'_s)=0$ would mean $\F'_s=\emptyset$, and consequently
$|\F|=|\F'|\le b(n,s)$, which would be a contradiction with our choice of $\F$. Hence,
$\nu(\F_s)\ge1$ and, so, for some $1\le q\le s-1$, we must have $\nu(\F'_q)=\nu(\F'_{q+1})$,
meaning that $\F'_q$ has property ONE. Let

$$q_0=\min\{q:\; \F'_q \mbox{ has } ONE\}.$$
We have $1\le q_0\le s-1$, but, in fact, $q_0$ is much smaller.
\begin{Fact}\label{q0}
 $q_0\le.2s+1$
\end{Fact}
\proof We first claim that
\begin{equation}\label{4}
n-q_0\le4(s-q_0+1)-1.
\end{equation}
Indeed, we have
$$n':=|V(\F_{q_0}')|=n-q_0\quad\mbox{and}\quad s':=\nu(\F_{q_0}')=s-q_0.$$ If (\ref{4}) would not
hold then $n'\ge4(s'+1)$, and by the result from \cite{FRR} mentioned in the Introduction, we would
have $|\F_{q_0}'|\le b(n-q_0,s-q_0)=\binom{n-q_0}3-\binom{n-s}3$ and thus,
$$|\F'|\le \binom n3-\binom{n-q_0}3+|\F_{q_0}'|\le \binom n3-\binom{n-s}3=M(n,s),$$
a contradiction. Consequently, (\ref{4}) holds, that is, $n-q_0\le 4(s-q_0)+3$, and since by Fact
\ref{n0s}(ii), $n\ge n_0(s)-1\ge 3.4s$, we have $q_0\le.2s+1$. \qed

\medskip

By Fact \ref{q0} and  (\ref{1}),
 $s'=s-q_0\ge.8s-1\ge s_0$.
Since $\F_{q_0}'$ has property ONE and $n'\ge 3s'+3$, by the assumption of Lemma \ref{L2},
$$|\F'_{q_0}|\le
m_{ONE}(n',s')\le M(n',s')= M(n-q_0,s-q_0).$$ We are going to show by inverse induction on $q$ that
for $q=q_0,\dots,0$, $|\F_q'|\le M(n-q,s-q)$. By estimate (iii) from Fact \ref{n0s}, for $q\ge1$,
$$n-q\ge n_1(s)-1-q\ge n_1(s-q)+2q-1-q\ge n_1(s-q)$$
and, consequently, $M(n-q,s-q)=b(n-q,s-q_)=\binom{n-q}3-\binom{n-s}3$. Thus, the inductive step, for  $1\le q\le q_0$, can be easily verified:
$$|\F_{q-1}'|\le \binom{n-q}2+|\F_q'|\le\binom{n-q}2+\binom{n-q}3-\binom{n-s}3=\binom{n-q+1}3-\binom{n-s}3.$$
The case $q=0$, that is, the inequality $|\F'_0|\le M(n,s)$, contradicts our assumption that
$|\F|=|\F'|=|\F_0'|>M(n.s)$. The proof of Lemma \ref{L2} is completed. \qed

\section{Proof of Lemma \ref{L1}}\label{P1}

Before we turn to the actual proof, we need to prove more facts about stable 3-graphs. Let $\F$ be a stable 3-graph with vertex set $[n]$, $n\ge 3s+3$.
Set $\nu(\F)=s$ and define
$$\F_0=\{F\cap[3s+2]:\; F\in\F\}.$$
Since, by stability, there is an $s$-matching in $\F$ with vertex set $[3s]$, no edge of $\F$ can
be disjoint from $[3s]$, and even more so from $[3s+2]$. Hence, $\emptyset\not\in\F_0$. Similarly,
if, in addition, $\F$ had property ONE, then there would be an $s$-matching in $\F$ with vertex set
$[2,3s+1]$ and so, no edge of $\F$ might share  just one vertex with $[3s+2]$. Hence,
\begin{equation}\label{H2}
\F\in ONE\quad\Rightarrow\quad \forall H\in\F_0:\;\; |H|\ge2.
\end{equation}

The following two observations play a crucial role in the proof of Lemma \ref{L1}.
\begin{Fact}\label{F0} If $\F$ is stable, then $\F_0$ is stable and
$\nu(\F_0)=s$.
\end{Fact}

\proof The stability of $\F_0$ follows directly from the stability of $\F$. Since $\F_0$ contains
an $s$-matching of $\F$, we also have $\nu(\F_0)\ge s$. If $\nu(\F_0)> s$, then there would exist
$s+1$ disjoint subsets $H_1,\dots,H_{s+1}$ of $[3s+2]$, each of size 1, 2 or 3, and $s+1$ triples
$F'_1,\dots,F'_{s+1}$ such that $H_i\subseteq F'_i$, $i=1,\dots,s+1$, and $F_i'\cap[3s+2]=H_i$. Let
us choose the $H_i$'s and  the $F'_i$'s  to maximize $|\bigcup_{i=1}^{s+1}F'_i|$. If the union had
fewer than $3s+3$ vertices, then some two triples, say $F'_1$ and $F'_2$, would intersect, and thus
there would exist a vertex $v\in [3s+2]\setminus\bigcup_{i=1}^{s+1}F'_i$. Then, denoting by $u$ a
common vertex of $F'_1$ and $F'_2$, we could replace $F'_1$ by $(F'_1\setminus\{u\})\cup\{v\}$,
obtaining a new family with a larger union, a contradiction.  Thus, $F_1',\dots, F'_{s+1}$ are
pairwise disjoint, which contradicts the assumption that $\nu(\F)=s$. \qed

\medskip

We say that $\F$ is \emph{maximal} if for every $E\not\in\F$,  $\nu(\F\cup\{E\})>\nu(\F)$.

\begin{Fact}\label{super} If $\F$ is stable and maximal, then $\F$ is closed under taking 3-element supersets of the sets of size two in $\F_0$.
\end{Fact}

\proof Let $H\in\F_0$, $|H|=2$. There exists $v\not\in [3s+2]$ such that $H\cup\{v\}\in\F$.  By
stability, also $H\cup\{v'\}\in\F$ for all $v'<v$, $v'\not\in H$. Suppose that there exists $u>v$
with $H\cup\{u\}\not\in \F$. By maximality, the only reason for that is that there is an
$(s+1)$-matching $F_1,\dots, F_{s+1}$ in $\F\cup\{H\cup\{u\}\}$. W.l.o.g. set $F_{s+1}=H\cup\{u\}$.

We are going to show that $H\cup\{u\}\not\in\F$ implies that one can replace $u$ by some $w$ with
$F_1,\dots, F_s, H\cup\{w\}$ forming a matching in $\F$, which is a contradiction. To find such $w$
observe that $|([3s+2]\setminus H)\cup\{v\}|=3s+1$, and hence there is $w\in ([3s+2]\setminus
H)\cup\{v\}$ not belonging to $\bigcup_{i=1}^sF_i$. But then, replacing $F_{s+1}$ with $H\cup\{w\}$
leads to an $(s+1)$-matching in $\F$. \qed

\subsection{Set-up}


Let integers $s$ and $n$ and a 3-graph $\F$ be such that

\begin{enumerate}
\item $s\ge25$,
\item $n\in\{n_1(s)-1,n_1(s)\}$,
\item $\nu(\F)=s$
\item $\F$ has property ONE (and so, $\F$ is stable),
\item $|\F|=m_{ONE}(n,s).$

\noindent In addition, as it is shown below, we may also assume that
\item $\F$ is maximal (with respect to $\nu(\F)=s$).
\end{enumerate}


 \begin{Fact}\label{max}
 If $\F$ satisfies (iii)-(v), then $\F$ is maximal.
 \end{Fact}
 \proof Observe that for each $E\not\in\F$,
$$\nu(\F\cup\{E\})\ge\nu(sh(\F\cup\{E\}))\ge\nu(\F)=s,$$
the second inequality due to the inclusion $\F\subset sh(\F\cup\{E\})$. Suppose that
\newline $\nu(\F\cup\{E\})=\nu(sh(\F\cup\{E\}))=s$. Then, $sh(\F\cup\{E\})$ has ONE, because $\F$ did.
But $|sh(\F\cup\{E\})|>|\F|=m_{ONE}(n,s)$ which is a contradiction. \qed

\bigskip

By Fact \ref{R}, to prove Lemma \ref{L1}, it is sufficient to show that a 3-graph $\F$ satisfying
(i)-(vi) above has at most $M(n,s)$ edges.
 As an immediate consequence of (vi) and Fact \ref{super}, we obtain a pivotal identity:
\begin{equation}\label{5}
|\F|=\sum_{H\in\F_0}\binom{n-3s-2}{3-|H|}.
\end{equation}
Recall (\ref{H2}) and let $\F_0^i$, $i=2,3$, stand for the subhypergraph of $\F_0$ consisting of
all edges of size $i$. Then, the above identity can be rewritten as
$$|\F|=|\F_0^3|+(n-3s-2)|\F_0^2|.$$
We are going to rewrite identity (\ref{5}) one more time. By property ONE, there is a matching of
size $s$ in $\F$ avoiding  vertex 1. Thus, by stability, there is also an $s$-matching contained in
$[3s+2]\setminus\{1\}$. Let $F_1,\dots,F_s$ form such a matching and let
$$F_0:=\{1,d\}=[3s+2]\setminus(F_1\cup\cdots\cup F_s)$$
contain the two remaining elements of $[3s+2]$. Among all possible candidates for $F_1,\dots,F_s$
we choose one which makes $d$ as small as possible.
  Observe that $F_0\not\in\F_0$, and consequently, for every $v\in[3s+2]$, $v\neq d$, we also have $\{d,v\}\not\in \F_0$.

 For $H\in\F_0$, we define two parameters: \emph{the spread}
$$z(H)=|\{i:\; H\cap F_i\neq\emptyset\}|$$
and \emph{the weight}
$$w(H)=\frac{\binom{n-3s-2}{3-|H|}}{\binom{s-z(H)}{3-z(H)}}.$$
For each triple of indices $\tau\in\binom{[s]}3$, let
$$V^\tau=F_0\cup\bigcup_{i\in\tau}F_i\quad\mbox{and}\quad\F^{\tau}_0=\{H\in\F_0:\; H\subset V^\tau\}.$$
It follows from (\ref{5}) and the definition of $w(H)$ that

\begin{equation}\label{6}
|\F|=\sum_{\tau\in\binom{[s]}3}\sum_{H\in\F^\tau_0}w(H).
\end{equation}

As $M(n,s)\ge a(s)=|\A|$, our ultimate goal is to show that $|\F|\le|\A|$.
 Recall that $\A$ is the complete 3-graph on
$[3s+2]$. Identity (\ref{6}) holds also for $\A$ instead of $\F$ (with the same choice of an
$s$-matching $F_1,\dots,F_s$), but due to the symmetry of $\A$, the inner sums in (\ref{6}) are
independent of the choice of $\tau$, and thus all equal to each other. Denoting this common value
by $W$, we thus have $|\A|=a(s)=\binom s3W$, and we will achieve our goal by showing that for each
$\tau\in\binom{[s]}3$

\begin{equation}\label{7}
\sum_{H\in\F^\tau_0}w(H)\le\sum_{H\in\A^\tau_0}w(H).
\end{equation}

\subsection{The 11-vertex board}

\begin{figure}[!ht]
\centering
\includegraphics [width=15cm]{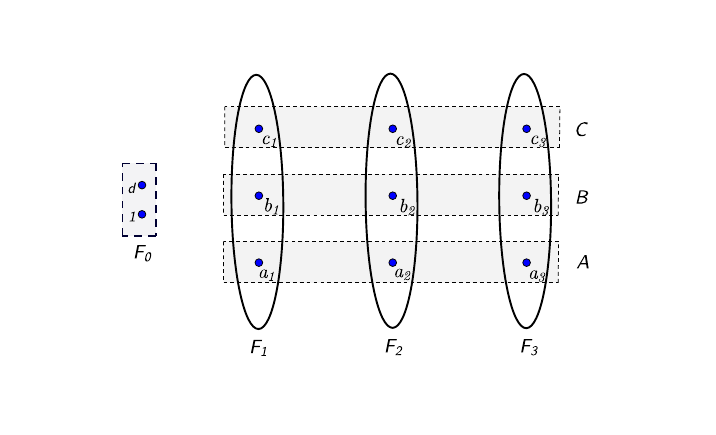}
\caption{The 11-vertex board $V^\tau$}\label{Fig1}
\end{figure}

Let us fix $\tau\in\binom{[s]}3$. Without loss of generality assume that $\tau=\{1,2,3\}$ and,
thus, $V^\tau=F_0\cup F_1\cup F_2\cup F_3$. Set $F_i=\{a_i<b_i<c_i\}$, $i=1,2,3$, and
$A=\{a_1,a_2,a_3\}, B=\{b_1,b_2,b_3\}, C=\{c_1,c_2,c_3\}$ (see Fig. \ref{Fig1}).


The family
$\F_0^\tau$ consists of triples and pairs. According to their width we divide the triples into
\emph{wide} ($z(H)=3$), \emph{medium} ($z(H)=2$), and \emph{narrow} ($z(H)=1$), with the
corresponding weights $w(H)$ equal to
$$1,\quad \frac1{s-2},\quad \mbox{and}\quad \frac1{\binom{s-1}2}.$$ Similarly, we divide the pairs into \emph{wide}
($z(H)=2$) and \emph{narrow} ($z(H)=1$), with weights
$$\frac{n-3s-2}{s-2} \quad \mbox{and}\quad \frac{n-3s-2}{\binom{s-1}2}.$$

Our strategy for proving (\ref{7}) is as follows. First note that for each triple $H\in\F_0$ its weight $w(H)$ is the same in $\F $ and in $\A$. Consequently, as $\A=\binom{[3s+2]}3$, for every triple $H\in\F_0^\tau$ its weights on both sides of (\ref{7}) cancel out. We will consider several cases with respect to the structure of $\F_0^\tau$, and in each case will argue that the weights $w(H)$ of the pairs  $H\in\F_0^\tau$ sum up to no more than the total weight of the triples which are absent from $\F_0^\tau$, therefore establishing (\ref{7}). Note that for a wide pair $H$, owing to the estimate $n\le n_1(s)\le 3.5s+3$ (see  Fact \ref{n0s}(i)), we can bound its weight by
$$w(H)=\frac{n-3s-2}{s-2}\le\frac{s/2+1}{s-2}.$$
Hence, for large $s$, it is enough to show that the number of wide pairs present in $\F_0^\tau$ is strictly less than twice the number of wide triples missing from $\F_0^\tau$. For smaller $s$, however, we need also look at  other types of sets (not just wide).

The absence of specified sets from $\F_0^\tau$ will be often forced  by the same kind of argument:
assuming their presence we would get a matching of size 4 in $\F_0^\tau$, leading to a
contradiction with Fact \ref{F0}. Not to repeat ourselves, we will refer to this argument as
\emph{a 4-matching argument}.

As a first example of this technique, consider a narrow pair $H$ which, say, is contained in
$F_0\cup F_1$. If its complement $H'=(F_0\cup F_1)\setminus H$, which is a narrow triple,
belonged to $\F_0^\tau$, then $H$, $H'$, $F_2$, and $F_3$ would form a 4-matching. Thus,
$H'\not\in\F_0^\tau$, and when proving (\ref{7}), we can use the bound
\begin{equation}\label{8}
w(H)-w(H')\le\frac{n-3s-3}{\binom{s-1}2}\le \frac{s}{(s-1)(s-2)},
\end{equation}
where the second inequality follows again from  Fact \ref{n0s}(i).

The next fact sets an upper bound on the total number of narrow pairs in $\F_0^\tau$.
\begin{Fact}\label{narrow}
For each $i\in\{1,2,3\}$, there are at most 3 narrow pairs in $F_0\cup F_i$.
\end{Fact}
\proof If $\{a_i,c_i\}\in\F_0^\tau$ then, by stability, also $\{1,b_i\}\in\F_0^\tau$, which is a contradiction by the 4-matching argument.
By stability again, this also excludes $\{b_i,c_i\}$ which majorizes $\{a_i,c_i\}$. Similarly, we cannot have both,  $\{1,c_i\}\in\F_0^\tau$ and  $\{a_i,b_i\}\in\F_0^\tau$. Finally, recall that no pair in $\F_0^\tau$ contains $d$. \qed

\medskip

Now we prove the existence of narrow triples  which often help to complete a 4-matching in $\F_0^\tau$.

\begin{Fact}\label{1db}
For each $i\in\{1,2,3\}$, we have $\{1,d,b_i\}\in\F_0^\tau$ (and thus, by stability, also $\{1,d,a_i\}\in\F_0^\tau$).
\end{Fact}
\proof If $b_i>d$, then $\{1,d,b_i\}\prec\{a_i,b_i,c_i\}=F_i\in\F_0^\tau$, and so, by stability,  $\{1,d,b_i\}\in\F_0^\tau$.
Assume now that $b_i<d$. We claim that $\{1,b_i\}\in \F_0$. Indeed, otherwise, by maximality of $\F$, there is in $\F$ an $s$-matching disjoint from $\{1,b_i\}$ which, by stability, implies the presence of such a matching in $[3s+2]\setminus\{1,b_i\}$. This contradicts the minimality of~$d$. Hence, $\{1,b_i\}\in \F_0$ and, by Fact \ref{super}, $\{1,b_i,d\}\in\F$. \qed

\subsection{The proof of (\ref{7})}

The proof of  (\ref{7}) is split into two major cases, quite uneven in length.

\medskip

\noindent{\bf Case I:} There exist $F',F''\in\F_0^\tau$ such that $F'\cup F''=B\cup C$.

\medskip

We claim there is no wide pair in $\F_0^\tau$. Indeed,  otherwise, by stability, there would be a pair $H\in\F_0^\tau$ with $H\subset A$. W.l.o.g.
assume that $a_1\not\in H$. Then $F',F'',H,$ and $\{1,d,a_1\}$ (by Fact \ref{1db}) would form a
4-matching in $\F_0^\tau$, a contradiction.
In view of the absence of wide pairs, our goal is to outweigh the contribution (to the L-H-S of (\ref{7})) of
the narrow pairs by the contribution (to the R-H-S of (\ref{7})) of the narrow and medium triples which
are missing from $\F_0^\tau$. The narrow triples are already taken care of via estimate (\ref{8}).
If there is a narrow pair in $\F_0^\tau$ which is contained in, say, $F_0\cup F_1$, then by
stability, also $\{1,a_1\}\in \F_0^\tau$. Consequently, a 4-matching argument, like the one used above, excludes from $\F_0^\tau$ the medium triple $\{d,a_2,a_3\}$ along
with all 8 of its majorants. Thus, (\ref{7})  holds if
$$\frac{3s}{(s-1)(s-2)}\le\frac9{s-2},$$ which is true for $s\ge3$.

\medskip

\noindent{\bf Case II:}  There are no sets $F',F''\in\F_0^\tau$ such that $F'\cup F''=B\cup C$.

\medskip

Note that if at least 5 of all 8 wide triples contained in $B\cup C$ were present in $\F_0^\tau$,
then Case I would hold. Thus,  there are at least 4 wide triples
missing from $\F_0^\tau$ (and contained in $B\cup C$). By Fact \ref{narrow} we have always at
most 9 narrow pairs and we will not attempt to refine this bound.
Instead, in each case we will try to bound from above the numbers of wide pairs  and,
from below, the number of wide triples, by, respectively, $x$ and $y$, so that
\begin{equation}\label{xyz}
x\times \frac{s/2+1}{s-2}+\frac{9s}{(s-1)(s-2)}\le y.
\end{equation}
In conclusion, if for some $s$, our bounds $x$ and $y$ satisfy (\ref{xyz}), then (\ref{7}) holds.

 Case II will have several subcases, for concise description of which we introduce the following notation:
$$BC=\{\mbox{wide pairs } H\in\F_0^\tau:\; H\cap B\neq\emptyset,\; H\cap C\neq\emptyset\},$$
$$\overline{BC}=BC\cup\left(\binom C2\cap\F_0^\tau\right),$$
and
$$\overline{\underline{BC}}=\overline{BC}\cup\left(\binom B2\cap\F_0^\tau\right).$$
Sets (of wide pairs) $AB, AC, \overline{AB}$, and $\underline{\overline{AB}}$ are defined analogously.

In each subcase below (except for the last one) we just derive the bounds $x$ and $y$ and claim that (\ref{xyz}) holds (for $s\ge25$), leaving numerical details to the untrusting reader.

\medskip

\noindent{\bf Subcase II.1:} $\overline{BC}\neq\emptyset$. W.l.o.g, let $\{c_1,b_2\}\in\F_0^\tau$. If $\{a_1,c_2,a_3\}\in\F_0^\tau$, then together with $\{c_1,b_2\}, \{b_1,a_2\}$, and $\{1,d,b_3\}$ it would form a 4-matching, a contradiction. Thus, $\{a_1,c_2,a_3\}\not\in\F_0^\tau$, there are at least $y=9$ wide triples (all majorants of $\{a_1,c_2,a_3\}$) missing from $\F_0^\tau$.

We claim
there are at most $x=12$ wide pairs in $\F_0^\tau$.
Indeed, by Fact \ref{super},
$$\{a_1,a_3\},\{a_1,c_2\},\{c_2,a_3\}\not\in\F_0^\tau,$$ which, by stability, excludes 9+3+3=15 out of all 27 wide pairs.
With $x=12$ and $y=9$, (\ref{xyz}) holds (in fact, already for $s\ge14$).

\medskip

\noindent{\bf Subcase II.2:} $\overline{BC}=\emptyset$. This subcase is further subdivided according to the size and structure of the graph $AC$.
Note that if $|AC|>2$ then $AC$ must have two disjoint edges.

\medskip

\noindent{\bf Subcase II.2.(i):} there is a 2-matching $M$ in $AC$. Say $a_1,a_2\in V(M)$. Then, by the 4-matching argument, $\{b_1,b_2,a_3\}\not\in\F_0^\tau$, and so, $y=12$. As we also have $x\le18$, (\ref{xyz}) holds. (By Fact \ref{super}, we can improve the upper bound on the number of wide pairs to $x=27-(4+6+6)=11$, but we do not need it.)

\medskip

\noindent{\bf Subcase II.2.(ii):} $|AC|=2$ and the two edges in $AC$ have a common endpoint in $C$, say $\{a_1,c_2\},\{a_3,c_2\}\in\F_0^\tau$.
Then, by stability, also $\{a_3,a_2\}\in\F_0^\tau$, which, by a 4-matching argument, excludes both $\{c_1,b_2,b_3\}$ and $\{b_1,b_2,c_3\}$ from $\F_0^\tau$.
This sets $y=6$. As for $x$, note that by a similar 4-matching argument, we may exclude from $\F_0^\tau$ both $\{b_1,a_2\}$ and $\{b_3,a_2\}$, and thus, by stability also $\{b_1,b_2\}$ and $\{b_2,b_3\}$. Hence, $|\underline{\overline{AB}}|\le8$ and so $x=10$ which yields (\ref{xyz}).

\medskip

\noindent{\bf Subcase II.2.(iii):} $|AC|=2$ and the two edges in $AC$ have a common endpoint in $A$, say $\{c_1,a_2\},\{c_3,a_2\}\in\F_0^\tau$.
We claim that
 $|\underline{\overline{AB}}|\le8$. Indeed, by a 4-matching argument, either $\{a_1,b_3\}$ or $\{b_1,a_3\}$ is not in $\F_0^\tau$. By symmetry, assume that  $\{a_1,b_3\}\not\in\F_0^\tau$. By a similar argument,  $\{a_1,b_2\}$ or $\{b_1,a_3\}$ is not in $\F_0^\tau$, and  $\{b_1,b_3\}$ or $\{a_1,a_3\}$ is not in $\F_0^\tau$. In addition, $\{b_1,b_3\}\not\in\F_0^\tau$, by stability. This yields at least 4 wide pairs missing from $\underline{\overline{AB}}$, and thus,
$|\underline{\overline{AB}}|\le12-4=8$.

If any of $\{b_1,c_2,a_3\}, \{a_1,c_2,b_3\}$, or $ \{c_1,b_2,a_3\}$ is missing from $\F_0^\tau$, then, by stability, at least 6 wide triples are not in $\F_0^\tau$. Thus, we have $x=10$, $y=6$, and (\ref{xyz}) follows.
If, on the other hand, all these three triples are present in $\F_0^\tau$, then a 4-matching argument and stability imply that $\overline{AB}=\emptyset$.
Thus, we have $x=5$ and $y=4$, and (\ref{xyz}) holds again.

\medskip

\noindent{\bf Subcase II.2.(iv):} $|AC|=1$, say $\{a_1,c_2\}\in\F_0^\tau$.
If either of $\{c_1,a_2,b_3\}$ and  $\{c_1,b_2,a_3\}$ is missing from $\F_0^\tau$, then $y=6$ and, by a 4-matching argument, at least 3 pairs from $\underline{\overline{AB}}$ are missing ($\{b_1,b_2\}$ or $\{a_2,b_3\}$, $\{b_1,b_3\}$ or $\{b_2,a_3\}$, and $\{b_2,b_3\}$ or $\{a_2,a_3\}$). Thus, $x=10$ and (\ref{xyz}) follows.
Otherwise, by a 4-matching argument, $\underline{\overline{AB}}\subseteq\binom A2\cup\{\{a_1,b_2\},\{a_1,b_3\}\}$, so $x=6$ and (\ref{xyz}) holds again.

\medskip

\noindent{\bf Subcase II.2.(v):} $AC=\emptyset$. In this case, the total number of wide pairs coincides with $|\underline{\overline{AB}}|$ and we consider several subcases with respect to this quantity.

\medskip

\noindent{\bf $|\underline{\overline{AB}}|\le6$.} In this case, $x=6$, $y=4$, and  (\ref{xyz}) holds.

\medskip

\noindent{\bf $|\underline{\overline{AB}}|=7$.} It can be easily checked by inspection that $\underline{\overline{AB}}$ contains a matching $M$ of size 2 such that $A\setminus V(M)\neq\emptyset$, say $a_1\not\in V(M)$. Then, by a 4-matching argument, $\{a_1,c_2,c_3\}\not\in\F_0^\tau$. Since, in addition, there are at least 4 wide triple in $B\cup C$ missing from $\F_0^\tau$, we have $y=5$, $x=7$, and (\ref{xyz}) holds again.

\medskip

\noindent{\bf $8\le|\underline{\overline{AB}}|\le 10$.} Eight edges in $\underline{\overline{AB}}$ guarantee two matchings of size 2 in $\underline{\overline{AB}}$, $M'$ and $M''$, and two distinct vertices $a',a''\in A$ such that $a'\not\in V(M')$ and $a''\not\in V(M'')$. Let $a'=a_1$ and $a''=a_2$. Then, by the 4-matching argument, $\{a_1,c_2,c_3\}\not\in \F_0^\tau$ and $\{c_1,a_2,c_3\}\not\in \F_0^\tau$. Hence, $x=10$, $y=6$, and (\ref{xyz}) follows.

\medskip

\noindent{\bf $|\underline{\overline{AB}}|=11$.} In this case  there are three matchings  of size 2 in $\underline{\overline{AB}}$, $M_1, M_2, M_3$, such that $a_i\not\in V(M_i)$, $i=1,2,3$, and (\ref{xyz}) follows with $x=11$ and $y=7$.

\medskip

\noindent{\bf $|\underline{\overline{AB}}|=12$.}  If we just repeated the argument from the previous subcase, we would  have $x=12$ and $y=7$, and (\ref{xyz}) would hold for $s\ge33$ only. To push it down to $s\ge25$, we need to refine our argument and turn for help to medium triples.
But this is easy: all 12 triples of the form $\{a_i,c_i,c_j\}$ or $\{b_i,c_i,c_j\}$, $i\neq j$, are forbidden in $\F_0^\tau$. Indeed, if, say, $\{a_1,c_1,c_2\}\in\F_0^\tau$, then we would get a 4-matching consisting of $\{a_1,c_1,c_2\}$, $\{a_2,b_3\}$, $\{b_2,a_3\}$, and $\{1,d,b_1\}$.
Thus, (\ref{7}) follows from
$$12\times \frac{s/2+1}{s-2}+\frac{9s}{(s-1)(s-2)}-\frac{12}{s-2}\le7,$$
which is true for $s\ge25$. The proof of Lemma 1, and thus the proof of Theorem \ref{main}, is completed.


\end{document}